\documentclass[final,leqno,onefignum,onetabnum]{siamltex1213}

\usepackage[dvips]{epsfig}

\usepackage{geometry}                
\usepackage[parfill]{parskip}    

\usepackage{amsmath, amssymb}

\usepackage{graphicx}
\usepackage{overpic}
\usepackage{overpic}
\usepackage{cancel}
\usepackage{rotating}
\usepackage{url}
\usepackage{color}
\usepackage[dvipsnames]{xcolor}

\usepackage{rotating}

\newcommand{\bA}{\mathbf{A}}
\newcommand{\bB}{\mathbf{B}}

\newcommand{\bg}{\mathbf{g}}
\newcommand{\bh}{\mathbf{h}}

\newcommand{\bw}{\mathbf{w}}

\newcommand{\bx}{\mathbf{x}}
\newcommand{\by}{\mathbf{y}}
\newcommand{\bu}{\mathbf{u}}
\newcommand{\bv}{\mathbf{v}}
\newcommand{\bz}{\mathbf{z}}
\newcommand{\bq}{\mathbf{q}}

\newcommand{\bK}{\mathbf{K}}

\definecolor{blue}{rgb}{0,0,1}
\definecolor{darkgreen}{rgb}{0,0.5,0}
\definecolor{red}{rgb}{1,0,0}
\definecolor{teal}{rgb}{0,0.5,0.7}

\graphicspath{{Figures/}}

\title{Generalizing Koopman theory to allow for inputs and control}

\author{Joshua L. Proctor\footnotemark[1]\thanks{Corresponding author (\email{joproctor@intven.com})}
\and Steven L. Brunton\footnotemark[2]\ \footnotemark[3]
\and J. Nathan Kutz\footnotemark[3]}

\begin{document}
\maketitle

\renewcommand{\thefootnote}{\fnsymbol{footnote}}

\footnotetext[1]{Institute for Disease Modeling Bellevue, WA 98004, United States}
\footnotetext[2]{Department of Mechanical Engineering, University of Washington, Seattle, WA 98195, United States}
\footnotetext[3]{Department of Applied Mathematics, University of Washington, Seattle, WA 98195, United States}

\begin{abstract}
We develop a new generalization of Koopman operator theory that incorporates the effects of inputs and control.  
Koopman spectral analysis is a theoretical tool for the analysis of nonlinear dynamical systems.  
Moreover, Koopman is intimately connected to Dynamic Mode Decomposition (DMD), a method that discovers spatial-temporal coherent modes from data, connects local-linear analysis to nonlinear operator theory, and importantly creates an equation-free architecture allowing investigation of complex systems. 
In actuated systems, standard Koopman analysis and DMD are incapable of producing input-output models; moreover, the dynamics and the modes will be corrupted by external forcing.
Our new theoretical developments extend Koopman operator theory to allow for systems with nonlinear input-output characteristics.  We show how this generalization is rigorously connected and generalizes a recent development called Dynamic Mode Decomposition with control (DMDc).  
We demonstrate this new theory on nonlinear dynamical systems, including a standard Susceptible-Infectious-Recovered model with relevance to the analysis of infectious disease data with mass vaccination (actuation).  
\end{abstract}




\section{Introduction}

We introduce a new method called Koopman with inputs and control (KIC) that generalizes Koopman spectral theory to allow for the analysis of complex, input-output systems.  Koopman operator theory, which is built on the seminal contribution of Bernard Koopman in 1931~\cite{Koopman:1931}, is a powerful and increasingly prominent theory that allows one to transform a nonlinear dynamical system into an infinite-dimensional, linear system~\cite{Koopman:1931,Mezic:2005,Rowley:2009}.  Linear operator theory~\cite{Friedman}, specifically eigenfunction expansion techniques, can then be used to construct solutions of the original system.  As such, Koopman theory is perhaps an early theoretical predecessor of what is now called {\em nonlinear manifold learning}, i.e. discovering nonlinear manifolds on which data live.  In Koopman theory, the data considered is generated from a nonlinear dynamical system and candidate manifolds are constructed from observables of the original state-space variables.  In our KIC innovation, we consider a nonlinear dynamical system with inputs and outputs, thus requiring a generalization of Koopman's original definition.  We demonstrate the method on a number of examples to demonstrate the effectiveness and success of the technique.
Importantly, the Koopman method is a data-driven, model-free method that is capable of constructing the best (in a least-square sense) underlying dynamics and control of a given system from data alone.  This makes it an attractive data-driven architectures in modern dynamical systems theory.

Although proposed more than eight decades ago, few results followed the original formulation by Koopman~\cite{Koopman:1931}.  This was partly due to the fact that there was no efficient way proposed to compute the Koopman operator itself.  Additionally, even if an algorithm had been proposed, there were no computers available to compute them in practice during that time period.  Interest was revived once again in 2004/5 by Igor Mezi\'c {\em et al}~\cite{Mezic:2004,Mezic:2005} who showed that Koopman theory could be used for the spectral analysis of nonlinear dynamical systems.   Two critical and enabling breakthroughs came shortly after.  In 2008/10, Schmid and Sessterhen~\cite{Schmid2008APS} and Schmid~\cite{Schmid:2010} proposed the Dynamic Mode Decomposition (DMD) algorithm for decomposing complex, spatio-temporal data, and in 2009, Rowley {\em et al.}~\cite{Rowley:2009} showed that the DMD was, in fact, a computation of the Koopman operator for linear observables.  Most recently, Tu {\em et al}~\cite{Tu:2014a} generalized and improved the DMD algorithm and definition to its current, state-of-the-art form.  The combined work of Mezi\'c, Rowley, Schmid and their co-workers thus laid the theoretical foundations that have led to the tremendous subsequent success of the DMD/Koopman method.  In a very short period of time since, DMD theory 
has been applied with great success to a broad set of domain sciences including complex fluid flows~\cite{Schmid:2009,Schmid:2010,Schmid:2011,Grilli:2012,Bagheri:2013,Tu:2014a,tissot2014model}, foreground/background separation in video streams~\cite{Grosek:2013}, epidemiology~\cite{Proctor:2015EP}, and neuroscience~\cite{Brunton2016}.  The theory also allows for critical enabling theoretical augmentations that can take advantage of compression and sparsity~\cite{Jovanovic2014sparsity,brunton:2014b,Gueniat2015pof}, multi-resolution/multi-scale phenomenon~\cite{Kutz2016siads}, de-noising~\cite{Dawson2014arxiv,Hemati2015arxiv}, data fusion~\cite{Williams2015epl}, extended and kernel DMD~\cite{Williams2014arxivA,Williams2015jnls}, and control~\cite{Proctor:2015DMDc}.  Indeed, our objective is to describe how Koopman operator theory can be generalized to include the analysis of input-output systems.  Further, we demonstrate how Koopman is fundamentally connected to Dynamic Mode Decomposition with control (DMDc), a recently developed extension of DMD for input-output systems \cite{Proctor:2015DMDc} which has already been successfully applied to model a rapidly pitching airfoil \cite{Dawson2015}.   

The rapid adoption of Koopman theory across a number of scientific and engineering fields~\cite{budivsic2012applied,Mezic:2013} is not surprising.   Its fundamental success stems from the fact that it is an {\it equation-free} method, relying on data alone to reconstruct a linear dynamical system characterizing the nonlinear system under consideration.  Such linear systems may be characterized using basic methods from ordinary differential equations and spectral analysis, as shown by Mezi\'c~\cite{Mezic:2005}.  The method can be applied to high-dimensional measurement data collected from complex systems where governing equations are not readily available; and the numerical instantiation of Koopman can be orders of magnitude faster than solving for solutions of PDEs with complex domains.  KIC inherits these advantageous characteristics, but extends the domain of applicability to input-output systems.     

The control of high-dimensional, nonlinear systems is a challenging task that is of paramount importance for applications such as flow control \cite{Brunton:2015Rev} and eradicating infectious diseases \cite{Proctor:2015EP}.  The construction of effective controllers typically rely on relatively few states, a computationally feasible model to implement, and fast solvers to minimize latencies introduced by computing estimates of the system \cite{AstromAuto2014}.   Further, control laws often rely on solving a single large Riccati equation $(\mathcal{H}_2)$ or iteratively through another set of equations $(\mathcal{H}_{\infty})$.  For modern engineering systems with high-dimensional measurement data and possibly high-dimensional input data, the requirements of the controllers are too restrictive.  Thus, most practical methods for handling these modern systems rely heavily on dimensionality-reduction techniques.  These {\it model reduction} techniques typically employ the singular value decomposition to discover low-dimensional subspaces where the dynamics evolve \cite{HLBR_turb}.  On these low-dimensional subspaces, controllers can be described, constructed, and implemented \cite{Moore:1981,ERA:1985,HLBR_turb,rowley:2004,rowley.05,RapisardaAuto2011,Winck2013,FLM:9380777}.  Further, this paradigm is exemplified in the classic method called balanced truncation which utilizes both the low-dimensional controllable {\it and} observable subspaces to produce a balanced, reduced-order model for control \cite{Moore:1981}.  Notably, balanced truncation has been extended and generalized to handle high-dimensional measurement data by a method called balanced proper orthogonal decomposition (BPOD), but the method requires a linear adjoint calculation~\cite{Lall:2002,willcox:2002,rowley.05,ilak.08}, which is not possible in many data-driven experiments.  

The models produced by BPOD have been previously demonstrated to be equivalent to the balanced input-output models produced by the Eigensystem Realization Algorithm (ERA), a method developed to be used on linear and low-dimensional systems \cite{Ma2011tcfd}.  ERA and the Observer Kalman Identification method (OKID) are apart of a class of methods developed for system identification \cite{ERA:1985,OKID:1991,ForgioneAuto2014}.  Similar to DMD and DMDc, system identification methods are inherently equation-free, acting only on measurement and input data.  In fact, the modal decomposition methods have been shown to be intimately connected to ERA, OKID, and other system identification methods called subspace identification methods such as Numerical algorithms for Subspace State Space System Identification (N4SID) \cite{Qin20061502,Tu:2014a,Proctor:2015DMDc}.  In this manuscript, we demonstrate how KIC reduces to DMDc for linear input-output systems.  KIC can be interpreted in terms of nonlinear system identification since the architecture allows for the analysis of nonlinear systems.  

The outline of the paper is as follows:  \S~\ref{s:back} describes the background on Koopman operator theory and its connections to DMD.  \S~\ref{s:KIC} describes the new development called KIC and the strong connections to DMDc.  The following section \S~\ref{sec:applications} presents a number of numerical examples including nonlinear input-output systems.  

\section{Background: Koopman and Dynamic Mode Decomposition}
\label{s:back}

Koopman operator theory and DMD are powerful and intimately connected methods for analyzing complex systems.  Data collected from numerical simulations, experiments, or historical records can be utilized by Koopman and DMD to extract important dynamic characteristics relevant for prediction, bifurcation analysis, and parameter optimization.  This section provides the mathematical background for Koopman operator theory, DMD, and how they are connected \cite{Mezic:2005,Schmid2008APS,Schmid:2009,Rowley:2009,Tu:2014a}.  

\subsection{The Koopman Operator for dynamical systems}

The {\it Koopman operator} is a linear operator defined for any nonlinear system \cite{Koopman:1931}.  Spectral analysis of this linear operator provides an analytic and numerical tool to analyze flows arising from nonlinear dynamical systems \cite{Mezic:2005,Rowley:2009,budivsic2012applied,Mezic:2013}.  In this section, we describe the background on Koopman operator theory.  

Consider the discrete nonlinear dynamical system:
\begin{align}
\bx_{k+1} = \mathbf{f}(\bx_k),
\label{eq.nonlinearf}
\end{align}
evolving on a smooth manifold $\mathcal{M}$ where $\bx_k \in \mathcal{M} $.  $\mathbf{f}$ is a map from $\mathcal{M}$ to itself, and $k$ is an integer index.  For most practical engineering problems, we consider our state and manifold to be $\bx \in \mathbb{R}^{n_x}$.  We could equivalently describe Koopman operator theory for continuous-time systems, but here we restrict to the discrete-time setting as most engineering problems collect discrete time data.  We also define a set of scalar valued observable functions $g: \mathcal{H} \rightarrow \mathbb{R}$, which forms an infinite-dimensional Hilbert space.  This space consists of the Lebesque square-integrable functions on $\mathcal{H}$.  The Koopman operator $\mathcal{K}$ acts on this set of observable functions: 
\begin{align}
\mathcal{K} g(\bx) \triangleq g(\mathbf{f}(\bx)).
\label{eq.operf}
\end{align}
The Koopman operator is linear and {\it infinite}-dimensional, as defined in Eq.~\eqref{eq.operf}.  The nonlinear dynamical system is often considered finite-dimensional, but can be infinite-dimensional.  The linear characteristics of the Koopman Operator allow us to perform an eigendecomposition of $\mathcal{K}$:
\begin{align}
\mathcal{K} \varphi_j(\bx) = \lambda_j \varphi_j(\bx),~~~ j = 1,2,\dots,\infty.  
\label{eq.evUf}
\end{align}
Consider a vector-valued observable function $\bg : \mathcal{M} \rightarrow \mathbb{R}^{n_y}$.  Using the infinite expansion shown in Eq.~\eqref{eq.evUf}, the observable $\bg$, and if the $n_y$ components of $\bg$ lie within the span of eigenfunctions $\varphi_j$, the observable can be rewritten: 
\begin{align}
\bg(\bx) = \left [ \begin{array}{c}g_1(\bx) \\ g_2(\bx)\\ g_3(\bx) \\ \vdots \\ g_{n_y}(\bx) \end{array} \right ]  = \sum_{j = 1}^\infty \varphi_j(\bx)\bv_j 
\label{eq.evUf2}
\end{align}
where the vector valued coefficients $\bv_j$ are called Koopman modes.  Measure-preserving flows, as original considered in \cite{Koopman:1931}, allow for a specific description of the Koopman modes based on projections of the observables on to the span of $\mathcal{K}$:
\begin{align}
\bg(\bx) =\sum_{j=1}^\infty \varphi_j(\bx)   \left [ \begin{array}{c} \langle \varphi_j,g_1 \rangle_{\mathcal{H}} \\  \langle \varphi_j, g_2 \rangle_{\mathcal{H}} \\ \langle \varphi_j,g_3 \rangle_{\mathcal{H}} \\ \vdots \\  \langle \varphi_j,g_{n_y} \rangle_{\mathcal{H}}   \end{array} \right ] = \sum_{j = 1}^\infty \varphi_j(\bx)\bv_j .
\label{eq.evUf2cc}
\end{align}
The Koopman operator $\mathcal{K}$ is defined for all observables functions.  We later denote a finite-dimensional approximation of the Koopman operator (from data) as $\mathbf{K}$.  Rearranging terms from Eqs.~(\ref{eq.operf}) and (\ref{eq.evUf}) provides a new representation of the observable function $g$ in terms of Koopman modes and the corresponding Koopman eigenvalues $\lambda_j$:
\begin{align}
\mathcal{K} \bg(\bx) = \bg(\mathbf{f}({\bx})) =  \sum_{j = 1}^\infty \lambda_j \varphi_j(\bx)\bv_j
\label{eq.evUf3}
\end{align}
where the Koopman eigenvalues provide the growth/decay and frequency content of each Koopman modes, $\bv_j$.  For DMD, $\varphi_j(\bx)$ is a constant and is typically absorbed in to each of the modes.  For a linear operator and identity observable functions, e.g. $\bg(\bx) = \bx$, the eigenfunctions $\varphi(\bx)$ can be shown to be the inner product of the state $\bx$ with the left eigenvectors of the linear Koopman operator $\bw_j$ \cite{Rowley:2009}.  

A significant amount of recent work has focused on the application of the correct observable functions $g$ in order to uncover a Koopman operator that describes the nonlinear vector field \cite{Williams2015jnls,Brunton:2015Koopm}.  In particular, expanding the measured state in to a set of augmented states that either capture nonlinearities, i.e. $\bx^2$, $\bx^3$, $\sin(\bx)$, etc, or using the eigenfunctions of the underlying system.  In the examples for this manuscript, we will utilize these ideas to explore KIC.   

\subsection{Koopman and DMD}
\label{ss:kdmd}
We describe how the Koopman operator theory connects to DMD, thus intimately connecting measurement data with Koopman spectral analysis. Here, we follow the recent description provided in \cite{Tu:2014a}.  We describe a set of internal states $\bx_k$ where $k = 1,2,\dots,m$ by their respective measurements provided by Eqs.~\eqref{eq.evUf2} and \eqref{eq.evUf3}:
\begin{align}
\mathbf{y}_k = \bg(\bx_k),~~~~~\bz_k = \bg(\mathbf{f}(\bx_k)).
\label{eq.data1}
\end{align}
Note that the set of states $\bx_k$ do not need to be from a single trajectory of the dynamical system \cite{Tu:2014a}.  Each of the measurements can be collected to form two large data matrices:
\begin{align}
\mathbf{Y} = \left[ \begin{array}{ccccc} | & | &  & | \\
\by_1 &\by_{2}& \dots & \by_{m} \\
 | & | &  & | \end{array} \right],~~
 \mathbf{Z}&= \left[ \begin{array}{cccc} | & |  &  & | \\
\bz_1 &\bz_{2}  & \dots & \bz_{m} \\
 | & |  &  & | \end{array} \right].
\label{eq.snapshots}
\end{align}
\begin{definition}{Dynamic Mode Decomposition: (Tu et al. 2014 \cite{Tu:2014a})}
The dynamic mode decomposition of the measurement pair $(\mathbf{Y},\mathbf{Z})$ is given by the eigendecomposition of $\mathbf{A}$ where $\mathbf{A}  \triangleq \mathbf{Z}\mathbf{Y}^\dagger$ and $\dagger$ is the pseudo inverse.
\end{definition}\\
\\
{\it Remark:}  The measurements $\by_1$, $\by_2$, $\dots$, $\by_m$ do not have to be sequentially sampled.  The important relationship is between the current measurement and the future measurement, for example $\by_1$ and $\bz_1$.  The states $\bx_i$ do not have to be from a single trajectory of $\mathbf{f}$, but can be from a sample of the phase space.  Of course, collecting data from an experiment or a historical records, often this data will be collected from a single trajectory.  \\
\\
We can then compute DMD modes from the measurement pair by finding eigenvectors and eigenvalues that satisfy the standard eigenvalue problem:
\begin{align}
\mathbf{A} \bv_j = \lambda_j \bv_j.
\label{eq.evProb}
\end{align}
Assuming the matrix $\bA$ has a full set of eigenvectors, each measurement column $\by_k$ can be represented by expanding by the eigenvectors of $\bA$:
\begin{align}
\bg(\bx_k) = \sum_{j = 1}^n c_{jk}\bv_j.
\label{eq.finiteexp}
\end{align}
If we have linearly consistent data, the relationship $\bA \by_k = \bz_k$ is satisfied allowing us to apply the operator $\bA$  to Eq.~\eqref{eq.finiteexp}:
\begin{subequations}
\begin{equation}
\bg(\mathbf{f}(\bx_k)) = \bz_k = \bA \sum_{j = 1}^n c_{jk}\bv_j, 
\end{equation}
\begin{equation}
~~~~~~~~~~~~~~~~~=  \sum_{j = 1}^n \bA c_{jk}\bv_j, 
\end{equation}
\begin{equation}
~~~~~~~~~~~~~~~~~=  \sum_{j = 1}^n \lambda_j c_{jk}\bv_j,
\end{equation}
\label{eq.finiteexp2}
\end{subequations}
In the case of linearly consistent data matrices, the DMD modes and eigenvalues of Eq.~\eqref{eq.finiteexp2} correspond to the Koopman modes of Eq.~\eqref{eq.evUf3}.  We refer the reader to \cite{Tu:2014a} for a more detailed description.  
\section{Generalizing Koopman to allow for inputs and control}
\label{s:KIC}
In this section, Koopman operator theory is generalized to allow exogenous inputs to the systems.  In the first subsection, we show how Koopman operator theory can be generalized to include inputs.  Then, we show how this formulation can be applied to linear systems.  We describe the connection of this analysis to the DMDc resulting in a perspective on KIC that describes how to define a different output space for the Koopman operator.  

\subsection{Koopman with inputs and control}
\label{ss:KIC}

\begin{figure}
\begin{center}
\begin{overpic}[width=1\textwidth]{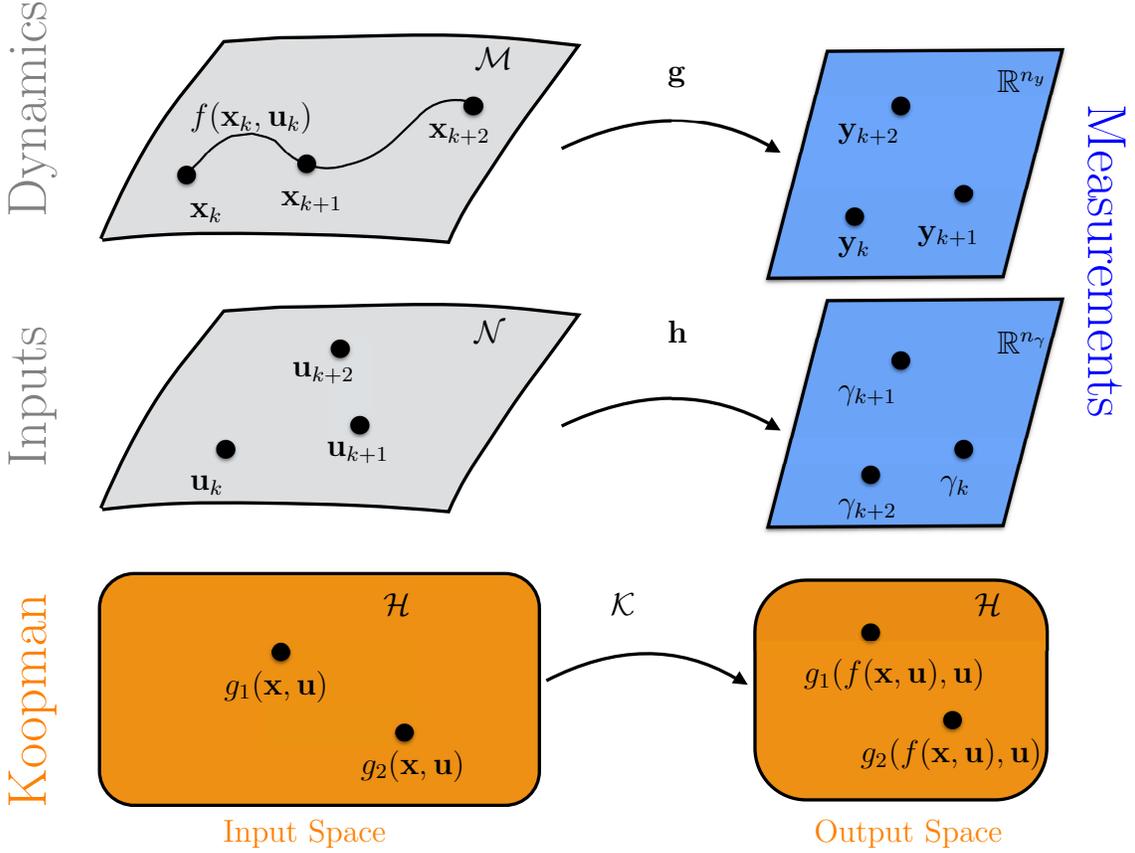}
	\large
	
	\put(92,65){\rotatebox{-90} {\huge \color{blue} Measurements}}
	\put(-3,55){\rotatebox{90} {\huge \color{gray} Dynamics}}
	\put(-3,33){\rotatebox{90} {\huge \color{gray} Inputs}}
	\put(-3,3){\rotatebox{90} {\huge \color{orange} Koopman}}
	\put(15,0){ {\large \color{orange} Input Space}}
	\put(67,0){ {\large \color{orange} Output Space}}

	\put(38,68){\large $\mathcal{M}$}	
	\put(38,44){\large $\mathcal{N}$}	
	\put(30,20){\large $\mathcal{H}$}
	
	\put(55,67){\large $\bg$}	
	\put(55,44){\large $\bh$}	
	\put(50,20){\large $\mathcal{K}$}
	
	\put(82,20){$\mathcal{H}$}	
	\put(84,43){$\mathbb{R}^{n_\gamma}$}	
	\put(84,66){$\mathbb{R}^{n_y}$}	

	\put(13,55){$\bx_k$}	
	\put(21,56){$\bx_{k+1}$}	
	\put(34,62){$\bx_{k+2}$}	
	
	\put(13,63){$f(\bx_k,\bu_k)$}	

	\put(13,31){$\bu_k$}	
	\put(25,34){$\bu_{k+1}$}	
	\put(22,41){$\bu_{k+2}$}	
	
	\put(70,52){$\by_k$}	
	\put(77,53){$\by_{k+1}$}	
	\put(70,62){$\by_{k+2}$}	
	
	\put(79,31){$\gamma_k$}	
	\put(70,39){$\gamma_{k+1}$}	
	\put(70,29){$\gamma_{k+2}$}	
	
	\put(16,13){$ g_1(\bx,\bu) $}
	\put(67,14){$ g_1(f(\bx,\bu),\bu) $}
	
	\put(28,6){$ g_2(\bx,\bu) $}
	\put(72,7){$ g_2(f(\bx,\bu),\bu) $}

	\end{overpic}
\end{center}
\caption{This figure illustrates the purpose of the Koopman operator with inputs and control.  The top row describes an underlying nonlinear dynamical system that is measured thought an observable function $\bg$.  The second row shows how the inputs, which can either be exogenous inputs or apart of a controller, are also measured.  The last row indicates the goal of the Koopman operator with inputs and control.  Namely, to find an operator that takes all observable functions $g_j(\bx,\bu)$ to the same observable function, but at a future internal state $g_j(f(\bx,\bu),\bu)$.    }\label{fig:illust1}
\end{figure}

Consider a nonlinear dynamical system that allows for external inputs
\begin{figure}
\begin{center}
\includegraphics[width=.9\textwidth]{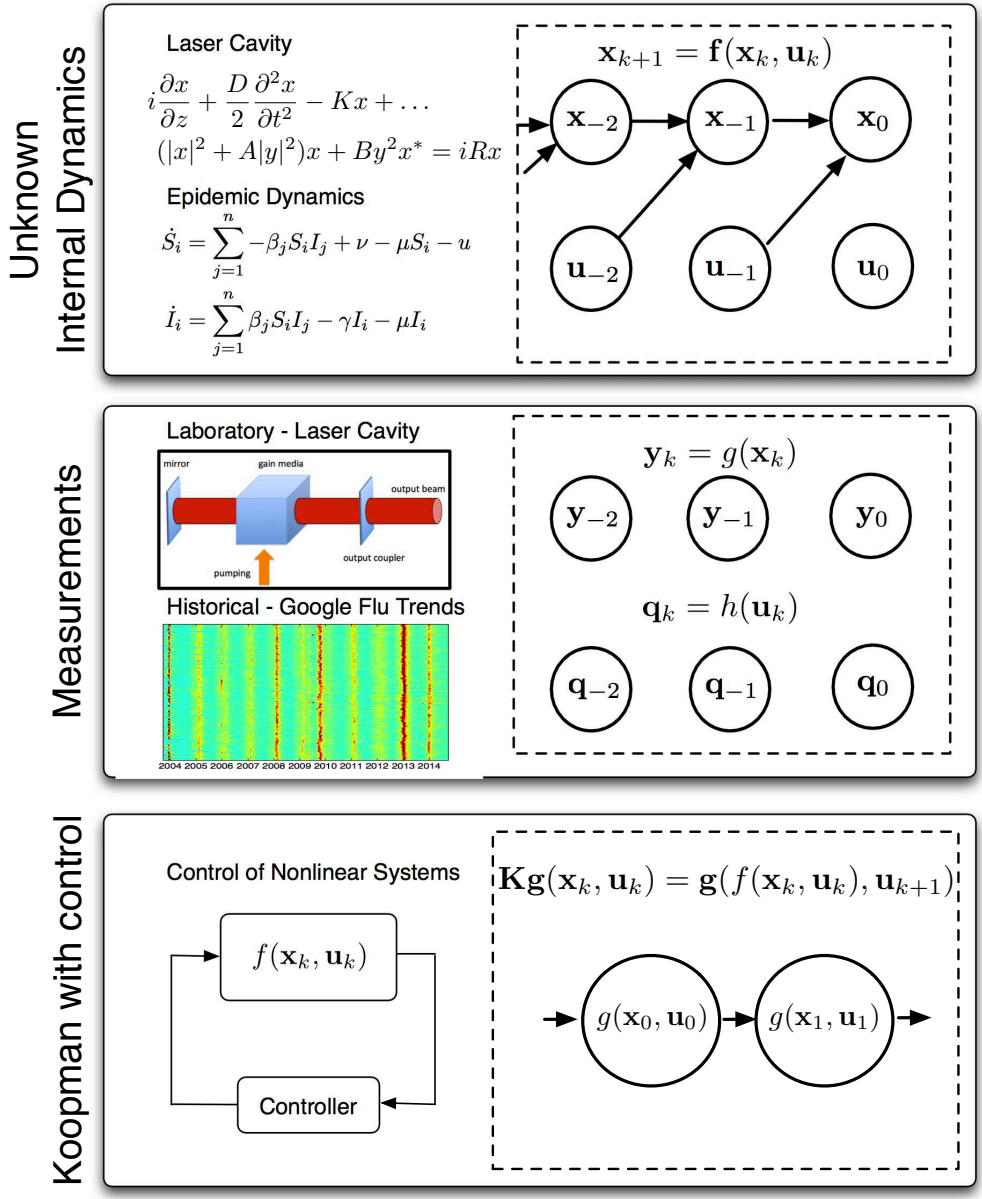}
\end{center}
\caption{An illustration about one of the goals of Koopman operator theory with or without inputs.  The first row shows that there might be an unknown system evolving according to some dynamical system.  The second row shows that we can measure the system experimentally, as in the case of optical systems, or historically, as in the case of historical infectious disease data.  The last row shows one of the goals of Koopman operator theory:  to discover an operator that can propagate forward in time a set of measurements for prediction and control.   }\label{fig:illust2}
\end{figure}
\begin{align}
\bx_{k+1} = \mathbf{f} (\bx_k,\bu_k),
\label{eq.nonlinearwc}
\end{align}
where $\bx \in \mathcal{M}$ and $\bu \in \mathcal{N}$ where both $\mathcal{M}$ and $\mathcal{N}$ are smooth manifolds.  As before, we dispense with the manifolds and consider $\bx \in  \mathbb{R}^{n_x}$ and $\bu \in \mathbb{R}^{n_u}$.  Further, we do not need $\bu$ to be constrained to a manifold.  We define a set of scalar-valued observation functions, but now the functions are dependent on the state {\em and} the input where $g: \mathcal{M} \otimes \mathcal{N} \rightarrow \mathbb{R}$.  Each observable function is an element of an infinite-dimensional Hilbert space $\mathcal{H}$.  Again, we choose the Hilbert space given by the Lebesque square-integrable functions.  Note that $\mathcal{H}$ can be broken in to three separate Hilbert spaces where the functions $g(\bx,\bu) = g(\bx)$ are in $\mathcal{H}_X$, $g(\bx,\bu) = g(\bu)$ are in $\mathcal{H}_U$, and finally the complement  $\mathcal{H}_{XU}$ which contain observable functions that offer mixed terms such as $g_{\bx,\bu} = x_1u_1$.  Thus, we can consider the Hilbert space to be composed of three components such as $\mathcal{H} = \mathcal{H}_X \otimes \mathcal{H}_U \otimes \mathcal{H}_{XU}$.  This partitioning could be extended to include linear identity observables, i.e. $g(\bx) = x_1$ where $x_1$ is the first element of $\bx$, versus nonlinear observables.   We take advantage of this construction later in this section to determine how the Koopman operator project to different partitions of $\mathcal{H}$ allowing us to connect KIC to DMDc.

The Koopman operator with inputs and control $\mathcal{K}: \mathcal{H}  \rightarrow \mathcal{H}$ acts on the Hilbert space of observable functions given by the following: 
\begin{align}
\mathcal{K} g(\bx,\bu) \triangleq g(\mathbf{f}(\bx,\bu),*).
\label{eq.defKwithcontrol}
\end{align}
where $*$ indicates a choice of definition.  Consider the following choices:
\begin{enumerate}
\item{$* = \bu$:  The inputs are evolving dynamically whether from state-dependent controllers or externally evolving systems such as those found in multi-scale modeling.   }
\item{$* =  \mathbf{0} $:  The inputs affect the state evolution, but the inputs are not evolving dynamically.  This is the case with impulse-response measurements and random exogenous disturbances.}
\end{enumerate}
The linear characteristics of the Koopman Operator allow us to perform an eigendecomposition of $\mathcal{K}$ given in the standard form:
\begin{align}
\mathcal{K} \varphi_j\left (\bx, \bu \right ) = \lambda_j \varphi_j\left (\bx, \bu \right ),~~~ j = 1,2,\dots.
\label{eq.expandKu}
\end{align}
The operator is now spanned by eigenfunctions that are defined by the inputs and state.  Using the infinite expansion shown in Eq.~\eqref{eq.evUf}, the observable functions $g_j$ can be rewritten in terms of the right eigenfunctions $\varphi_j$,
\begin{align}
 \bg(\bx,\bu)=
 \left [ \begin{array}{c}g_1(\bx,\bu) \\ g_2(\bx,\bu) \\ \vdots \\ g_{n_y}(\bx,\bu) \end{array} \right ]   =
 \sum_{j = 1}^\infty \varphi_j \left (  \bx,  \bu \right) \bv_j,
\label{eq.defineObservMM}
\end{align}
where $n_y$ is the number of measurements.  The new Koopman operator can be applied to this representation of the measurement
\begin{align}
\mathcal{K}  \bg(\bx,\bu)  = \bg(\mathbf{f}(\bx,\bu),\bu) =  \sum_{j = 1}^\infty \lambda_j  \varphi_j \left (  \bx ,\bu  \right)  \bv_j.
\label{eq.defineObservWithK}
\end{align}
Note that the expansion is in terms of Koopman eigenfunctions with vector valued coefficients that we call Koopman modes $\bv_j$.   The terminology of Koopman operator theory now allows for measurement functions that accept inputs.\\
\\
{\it Remark 1:}  In (\ref{eq.defKwithcontrol}), we could use $*=\mathbf{0}$ for the definition if we are not attempting to discover dynamics for the inputs.  Considering the discrete dynamical system of (\ref{eq.nonlinearwc}), the definition could also be $\mathcal{K} g(\bx_k,\bu_k) \triangleq g(\mathbf{f}(\bx_k,\bu_k),\bu_k)$ where the operator $\mathcal{K}$ will discover an identity map from $\bu_k$ to $\bu_k$ instead of a map from $\bu_k$ to $\mathbf{0}$.  This choice requires some a-priori information about the system and helps define a {\it family} of Koopman operators for the set of observable functions. This general perspective is discussed in more detail in \S\ref{ss:kio}.  

{\it Remark 2:}  In (\ref{eq.defKwithcontrol}), we could use $*=\bu$ for the definition if there is prior information that the inputs are evolving according to a set of dynamics.  The discrete dynamical system of (\ref{eq.nonlinearwc}) would define the following operator:  $\mathcal{K} g(\bx_k,\bu_k) \triangleq g(\mathbf{f}(\bx_k,\bu_k),\bu_{k+1})$. In this case, $\bu_k$ could technically be adjoined to the state $\bx_k$ and the original definition of Koopman applied.  We believe that partitioning $\mathcal{H}$ according to state and input observables also helps partition the operator $\mathcal{K}$ to disambiguate the impact of the state dynamics and the inputs.  Further, if the system being measured is multi-scale in nature such that the inputs to one scale of the system are evolving according to their own dynamics, then we could define the operator as $*=h(\bu)$.  The Koopman operator would then discover how $\bu$ is evolving on one scale as well as how $\bx$ is evolving with $\bu$ as an input. 

\subsection{KIC for linear systems }\label{unlinearsys}
In this subsection, we demonstrate how this new definition of KIC with $*=\bu$ can be applied to linear systems.  Consider the linear dynamical system with inputs
\begin{align}
\bx_{k+1} = \bA \bx_k + \bB \bu_k.
\label{eq.linearEx}
\end{align}
We consider full-state access and full-input access by choosing observable functions that are the identity i.e. $\bg(\bx) = \bx$.  The linear dynamical system can be rewritten in terms of a new state $\bz$
\begin{subequations}
\begin{align}
~~~~~~~~~~~~~ \left [ \begin{array}{c}   \bx_{k+1} \\  \bu_{k+1} \end{array} \right ]  &=  \left [ \begin{array}{cc} \mathbf{G}_{11} &  \mathbf{G}_{12}\\  \mathbf{G}_{21}  & \mathbf{G}_{22} \end{array} \right ]   \left [ \begin{array}{c} \bx_k \\  
 \bu_k \end{array} \right ] ,
  \end{align}
 \begin{align}
 \bz_{k+1 }&= ~~\mathbf{G} ~~\bz_k .
 \end{align}
\label{eq.linearEx2}

\end{subequations}
The eigenvalues of $\mathbf{G}$ are also the eigenvalues of $\mathcal{K}$ and the left and right eigenvectors of $\mathbf{G}$ are related to the eigenfunctions of $\mathcal{K}$.  The description of this linear system for an input-output system is clearly not canonical.   Typically, the future state would not include the future input.  There are exceptions, though, especially when considering a common method for analyzing non-autonomous dynamical systems where time is treated as a state (creating an augmented state) and a vector field $f$ augmented with a simple ODE, $\dot{t}=1$ \cite{gucken.02}.   Further, the inputs may actually have dynamics.  Later in this section, we comment on modifying the Koopman operator with inputs and control to illustrate a more canonical view of an input-output system, thus connecting previous work on DMDc \cite{Proctor:2015DMDc}.   The decomposition of $\mathbf{G}$, with eigenvalues $\lambda_j$ and eigenvectors $\bv_j$ is 
\begin{align}
 \mathbf{G} \mathbf{v}_j= \lambda_j \bv_j,~~~~~j=1,2,\dots,n
\label{eq.SVDexpansion}
\end{align}
The state can be represented by an expansion in terms of the singular vectors $\bv$:
\begin{align}
 \bz = \sum_{j = 1}^n  \varphi_j \left ( \bx,\bu \right) \bv_j  =  \sum_{j = 1}^n  \langle {\bz,\bw_j} \rangle_{\mathcal{H}}  \bv_j
\label{eq.defineObservKu2}
\end{align}
where we specify the Hilbert space $\mathcal{H}$ that the inner product is defined on $\langle \cdot \rangle_\mathcal{H}$.  Also, $\bw_j$ are the left eigenvectors of the operator $\mathbf{G}$. Further, the eigenfunctions $\varphi_j$ are projections of the state on the eigenvectors $\bw$.  For linear systems, the Koopman operator is equivalent to the linear map $\mathbf{G}$.  Further, the Koopman modes (both the left and right) coincide with the eigenvectors of $\mathbf{G}$.
\\
\\
{\em Remark 1:}  The analysis of the linear system with the new definition of the Koopman operator illustrates how the previous definition of the Koopman operator can be extended to handle inputs.  However, the example also indicates a challenge presented by this choice of methodology.  We may not be interested in finding a Koopman operator that predicts the future input $\bu_{k+1}$.  Further, if the inputs are random disturbances or exogenous inputs, there will not likely be an operator that can predict the future input.  Thus, $\mathbf{G}_{21}$ and $\mathbf{G}_{22}$ create issues for solving for the approximate Koopman operator.  There are a number of ways to alleviate these challenges, which are addressed in \S\ref{ss:kio}.  \\

\subsection{KIC and connections to DMDc}\label{DMDc}

In this subsection, we describe how KIC is connected to DMDc.  As stated in the previous subsection, the definition of KIC does not appear to fit with the canonical view of linear input-output systems.  In this subsection, we will illustrate the flexibility of the new definition by demonstrating how to connect the new theory with a recently developed method called DMDc \cite{Proctor:2015DMDc}.  This connection parallels the link between Koopman operator theory and DMD \cite{Rowley:2009}. 

Similar to \S\ref{ss:kdmd}, we describe a set of internal states $\bx_k$ where $k = 1,2,\dots,m$ and now with a set of internal inputs $\bu_k$ with linear, identity measurements given by the following:
\begin{align}
\left [ \begin{array}{c} \mathbf{y}_k \\ \mathbf{\gamma}_k \end{array} \right ] =  \bg(\bx_k,\bu_k),~~~~~\left [ \begin{array}{c} \mathbf{z}_k \\ \mathbf{\delta}_k \end{array} \right ] = \bg(\mathbf{f}(\bx_k,\bu_k),\bu_{k+1}).
\label{eq.data2}
\end{align}
As with Exact DMD, the set of states $\bx_k$ and inputs $\bu_k$ do not need to be from a single trajectory of the dynamical system \cite{Tu:2014a}.  Each of the measurements can be collected to form two large data matrices, described by the following:
\begin{align}
\mathbf{\Omega}&=\left[ \begin{array}{c} \mathbf{Y} \\ \mathbf{\Upsilon} \end{array}\right ] 
= \left[ \begin{array}{cccc} | & |  &  & | \\
\by_1 &\by_{2}  & \dots & \by_{m} \\
 | & |  &  & | \\
 | & | &  & | \\
 \mathbf{\gamma}_1 &  \mathbf{\gamma}_2 & \dots &  \mathbf{\gamma}_m \\
 | & | & & |  \end{array} \right],~~~~~
  \mathbf{\Delta} = \left[ \begin{array}{c} \mathbf{Z} \\ \mathbf{\Xi} \end{array}\right ] =  \left[ \begin{array}{cccc} | & |  &  & | \\
\bz_1 &\bz_{2}  & \dots & \bz_{m} \\
 | & |  &  & | \\
 | & | &  & | \\
 \mathbf{\delta}_1 &  \mathbf{\delta}_2 & \dots &  \mathbf{\delta}_m \\
 | & | & & |  \end{array} \right].
 \label{eq.snapshots2}
\end{align}
\begin{definition}{Dynamic Mode Decomposition with control : (Proctor et al. 2015 \cite{Proctor:2015DMDc})}  The dynamic mode decomposition of the measurement trio $(\mathbf{Z},\mathbf{Y},\mathbf{\Upsilon})$ is given by the eigendecomposition of the operator $\bA$ where $\mathbf{\tilde{G}} = [\bA~~~ \bB]$ and $\mathbf{\tilde{G}}  \triangleq \mathbf{Z}\mathbf{\Omega}^\dagger$.  
\label{def:dmdc}
\end{definition}\\
\\
The DMDc is defined for three measurement matrices $(\mathbf{Z},\mathbf{Y},\mathbf{\Upsilon})$ providing a non-square operator $\mathbf{\tilde{G}}$ that helps identify input-output characteristics.  The DMD modes from the measurement trio can be found by performing the singular value decomposition:
\begin{align}
\mathbf{\tilde{G}} \bv_j = \sigma_j \bq_j.
\label{eq.evProbt2}
\end{align}
Assuming the matrix $\mathbf{\tilde{G}}$ has a full set of singular vectors, each measurement column of $\mathbf{\Omega}$ can be represented by expanding by the eigenvectors $\bv_j$:
\begin{align}
\bg(\bx_k,\bu_k) = \left [ \begin{array}{c} \by_k \\ \mathbf{\gamma}_k \end{array} \right ]=  \sum_{j = 1}^n \varphi_j \bv_j.
\label{eq.finiteexpt}
\end{align}
If we have linearly consistent data, the relationship $\mathbf{\tilde{G}} \left [ \begin{array}{c} \by_k \\ \mathbf{\gamma}_k \end{array} \right ]= \bz_k$ is satisfied allowing us to apply the operator $\mathbf{\tilde{G}}$  to Eq.~\eqref{eq.finiteexp} giving:
\begin{subequations}
\begin{align}
 \bz_k &= \mathbf{\tilde{G}} \sum_{j = 1}^n c_{jk}  \bv_j, 
 \end{align}
 \begin{align}
~~~~=  \sum_{j = 1}^n \sigma_j c_{jk} \bq_j,
\end{align} 
\label{eq.finiteexpt2}
\end{subequations}
Note the difference between DMDc and KIC.  In DMDc, the expansion can be in terms of either the input or output space since the operator $\tilde{G}$ is not square.  In the previous subsection, we illustrated how $\mathbf{G}$ is effectively square which can step forward not only observables on the state, but also for the inputs.   In the following section, we describe how we synthesize these two perspectives.

\subsection{Adapting KIC to allow for different input and output spaces}
\label{ss:kio}

In this subsection, we connect the KIC architecture to DMDc.  Further, we illustrate how the Koopman operator can be viewed as projecting from the complete Hilbert space $\mathcal{H}$ to a subspace of $\mathcal{H}$.  This perspective of the Koopman operator with inputs allows for a different input and output space, thus facilitating the connection to not only DMDc, but also to recent developments such as Kernel DMD \cite{Williams2014arxivA,Williams2015jnls} and Sparse Identification of Nonlinear Dynamics (SINDy) \cite{Brunton:2015Sparse}.  We begin with slightly different definitions of the Koopman operator itself to show the connections to DMDc.  We then demonstrate how the Koopman operator can be viewed in terms of different input and output spaces.  

\subsubsection{The inputs are not dynamically evolving}
\label{sss:No}

The Koopman operator in (\ref{eq.defKwithcontrol}) with $*=\mathbf{0}$ is
\begin{align}
\mathcal{K} g(\bx,\bu) \triangleq g(\mathbf{f}(\bx,\bu),0).
\label{eq.defKwithcontrolMod}
\end{align}
In this case, the operator is no longer attempting to fit a future input prediction. Instead, this new modified Koopman operator is only attempting to propagate the observable functions at the current state and inputs to the future observable functions on the state.  Numerically, this definition modifies (\ref{eq.linearEx2}) for linear systems so that
\begin{align}
 \left [ \begin{array}{c}   \bx_{k+1} \\  \bu_{k+1} \end{array} \right ]  &=  \left [ \begin{array}{cc} \mathbf{G}_{11} &  \mathbf{G}_{12} \\ 0 & 0 \end{array} \right ]   \left [ \begin{array}{c} \bx_k \\  \bu_k \end{array} \right ] , 
 \end{align}
 which can be reduced to 
  \begin{subequations}
  \begin{align}
 ~~~~~~~~~~~~~~\left [ \begin{array}{c}   \bx_{k+1} \end{array} \right ]  &=  \left [ \begin{array}{cc} \mathbf{G}_{11} &  \mathbf{G}_{12} \end{array} \right ]   \left [ \begin{array}{c} \bx_k \\  \bu_k \end{array} \right ] , 
 \end{align}
 \begin{align}
 \bx_{k+1 }&= ~~\mathbf{\tilde{G}} ~~\bz_k .
 \end{align}
\label{eq.linearEx2Mod2}
\end{subequations}
This interpretation of the Koopman operator connects the non-canonical form found in \S\ref{unlinearsys} with the canonical version of system identification methods described in \S\ref{DMDc}.  This construction of the Koopman operator forces a closer inspection of the eigenfunction expansion in (\ref{eq.expandKu}).  There is no longer a requirement for having equivalent eigenfunctions $\varphi_j$ for both the input and output spaces of the operator $\mathcal{K}$.  Here, the eigenfunctions $\varphi_j$ could be mapped to a restricted subspace of $\mathcal{H}$ that only concerns the prediction of the future state (without the future input). 

\subsubsection{Input and output spaces for the Koopman operator}
\label{ss:KICgen}

We investigate how the output space of the Koopman operator can be restricted to a subspace of $\mathcal{H}$.  In \S\ref{ss:KIC}, we illustrated how the Koopman operator is defined on $\mathcal{H}$ for all observable functions $g(\bx,\bu)$.  This space $\mathcal{H}$ can be broken up in to different subspaces.  We illustrate how these subspaces can be utilized to describe the output space of the Koopman operator.  We could expand the input and output spaces of $\mathcal{K}$ by the following:
\begin{align}
\mathcal{K} \varphi_j\left ( \bx,\bu \right ) = \sigma_j \psi_j(\bx,\bu),~~~ j = 1,2,\dots.
\label{eq.expand1}
\end{align}
where $\psi_j$ are eigenfunctions that span a subspace $\mathcal{H}_X$.  These are a subset of the eigenfunctions of $\mathcal{H}$.  The span of $\mathcal{H}$ includes the eigenfunctions of the linear measurements on $\bx$, but also the nonlinear measurements on $\bx$ and $\bu$.   The vector of observable functions $\bg(\bx,\bu)$ can still be defined as in (\ref{eq.defineObservMM}), but now the Koopman operator applied to this set of observables becomes 
\begin{align}
\mathcal{K} \bg(\bx,\bu) =  \sum_{j = 1}^n \mathcal{K} \varphi_j \left ( \bx,\bu \right) \left [ \begin{array}{c}q_1 \\ q_2 \\ \vdots \\  q_{nx} \\ 0 \\ \vdots \\ 0
 \end{array} \right ] \approx  \sum_{j = 1}^\infty \sigma_j\psi_j \left (  \bx \right)   \bq_j 
\label{eq.defineObservKu1}
\end{align}
where $n_x$ is a smaller set of observable functions than $n_y$ and $\bq_j$ are left Koopman modes.  This allows for their to be different input and output spaces for the Koopman operator expansion.  The distinction allows for the practitioner to investigate how the Koopman operator projects a large input space of observables that includes linear, nonlinear, and mixed terms to a restricted output space of only linear observables.  Clearly in the case of DMDc, measurements of the state and input are collected, but the Koopman analog is only focused on determining the future measurements on the state (with the impact of the input included).  This perspective expands this view to allow for many more types of measurements than in the DMDc context.  Namely, the input space can include a large functional library expanding the measurement set as in the work by Williams et al. \cite{Williams2015jnls} which has been shown to discover the underlying nonlinear dynamics with more clarity.  

In this framework the input space of the Koopman operator can be expanded in linear measurements of the state and input, nonlinear measurements of the state and inputs,  and mixed state-input terms.  The output space, though, can be restricted to a set that spans, for example, the linear state measurements.   \\
\\
{\bf On the linear case for different input / output spaces of the Koopman operator}\\
\\
We are only interested in how the Koopman operator maps to observables that are of the form $g(\bx,\bu) = g(\bx)$.  This is exactly the case for DMDc.  For observables that are the identity on the state and input measurements, the input space can then be represented by a similar expansion of \S\ref{ss:KIC}, in terms of the right Koopman modes $\bv_j$:
\begin{align}
 \left [ \begin{array}{c} \bx \\ \bu \end{array} \right ] = \bz=  \sum_{j = 1}^n  \varphi_j \left (\bz \right) \bv_j  =  \sum_{j = 1}^n  \langle {\bz,\mathbf{v}_j} \rangle_{\mathcal{H}} \bv_j
\end{align}
The Koopman operator with control $\mathcal{K}$ can be applied to the input state $\bz$.   
\begin{subequations}  
\begin{align}
 \mathcal{K} \bz &= \sum_{j = 1}^n  \langle {\mathcal{K} \bz,\bq_j} \rangle_{\mathcal{H}_X} ~\bq_j 
 \end{align}
 \begin{align}
 &= \sum_{j = 1}^n  \langle {\bz,\mathcal{K}^*\bq_j} \rangle_{\mathcal{H}}~  \bq_j 
  \end{align}
 \begin{align}
 &= \sum_{j = 1}^n  \langle {\bz, \sigma_j^* \bv_j} \rangle_{\mathcal{H}} ~ \bq_j 
  \end{align}
 \begin{align}
 &= \sum_{j = 1}^n  \sigma_j \langle {\bz, \bv_j} \rangle_{\mathcal{H} }~  \bq_j
\end{align}
\label{eq.DeriveLinear}
\end{subequations}
where now the output space is expanded by $\bq_j$. The Koopman operator and the DMDc operator are equivalent.  This required restricting the output space to a subspace $\mathcal{H}_X$ of $\mathcal{H}$.  

\section{Applications}\label{sec:applications}
This section explores the theoretical development of KIC on various linear and nonlinear examples. For examples 1-3, we assume the perspective of the applied scientist which will be taking a finite set of measurements that are intuitive for their system.  The first example shows KIC when we assume there are dynamics on the inputs.  The second example explore a nonlinear dynamical system with a quadratic nonlinearity well-studied in the Koopman and DMD literature where we assume there is not dynamics on the inputs.  The final example looks at a more difficult example inspired by the study of infectious disease.  This example illustrates a difficulty facing the community applying Koopman and extended DMD on certain types of nonlinear problems.  It also illustrates a potential advantage of this framework.
\subsection{Example 1 -- Linear system with inputs}\label{sec:results:linearsys}
Consider the following linear dynamical system:
\begin{align}
\left [ \begin{array}{c} x_1 \\ x_2 \end{array} \right ]_{k+1} = 
 \left [ \begin{array}{c} \mu x_1 \\ \lambda x_2 + \delta u \end{array} \right ]_k .
\end{align}
A similar example can be found in \cite{Proctor:2015DMDc}. If $|\lambda|$ and/or $|\mu|$ is $> 1$, the system is unstable.  The goal is to recover the underlying dynamics and input matrix when there are various types of inputs including random disturbances, a state-feedback controller, or a multi-scale system.  We assume full access to the state and inputs giving the following relationship between the states, inputs, and  measurements:
\begin{align}
\left [ \begin{array}{c} y_1 \\ y_2 \\ \gamma  \end{array} \right ] = 
\left [ \begin{array}{c} x_1 \\ x_2  \\ u \end{array} \right ] ,~~~~~~
\left [ \begin{array}{c} y_1 \\ y_2 \\  \end{array} \right ]_{k+1} =
 \left [ \begin{array}{cc} \mu & 0  \\ 0 & \lambda  \end{array} \right ] \left [ \begin{array}{c} y_1 \\ y_2  \end{array} \right ]_k + \left [ \begin{array}{c} 0 \\ \delta \end{array} \right ] \gamma_k
 \label{eq.gv}
\end{align}
The dynamical system can be rewritten in the KIC form with the definition $* = \bu$ where we are interested in finding dynamics for the inputs:  
\begin{align}
\left [ \begin{array}{c} y_1 \\ y_2 \\ \gamma \end{array} \right ]_{k+1} =
 \left [ \begin{array}{ccc} \mu & 0  & 0 \\ 0 & \lambda & \delta \\ a & b & c \end{array} \right ] \left [ \begin{array}{c} y_1 \\ y_2  \\ \gamma \end{array} \right ]_k .
\end{align}
where $a$, $b$, and $c$ depend on the types of inputs.  We first investigate when the inputs are random disturbances.  We collect measurements of the state and inputs to investigate the reconstruction of a finite-dimensional Koopman operator.  The following is the first five snapshots of a single realization:
\begin{subequations}
 \begin{align}
 \mathbf{\Omega} &= \left [ \begin{array}{ccccc} 5 & 0.5  & 0.05 & 0.005 & 0.0005 \\
2 & 2.999  & 4.497 & 6.749 & 10.132 \\
-0.001 & -0.001 & 0.002 & 0.009 & 0.004
\end{array} \right ] 
\end{align}
\begin{align}
\mathbf{\Delta} &= \left [ \begin{array}{ccccc}  0.5  & 0.05 & 0.005 & 0.0005 & 0.00005 \\
2.999  & 4.497 & 6.749 & 10.132 & 15.203 \\
-0.001 & 0.002 & 0.009 & 0.004 & 0.006
\end{array} \right ] 
\end{align}
\end{subequations}
The parameters used for this example are $\mu = 0.1$, $\lambda = 1.5$, and $\delta = 1$, where the linear system is unstable.   The random disturbances for the input are zero mean and gaussian distributed with a variance of $0.01$  Six snapshots of data are used for the computation.  Using these data matrices, a restricted Koopman operator can be constructed, see (\ref{eq.linearEx2}) for an example.  The solution using these data matrices is:
\begin{align}
\mathbf{G} &=\left[ \begin{array}{cc} \mathbf{G}_{11}& \mathbf{G}_{12} \\  \mathbf{G}_{21}& \mathbf{G}_{22}  \end{array}\right] \approx
 \left [ \begin{array}{cc} \begin{bmatrix}0.1  & 0  \\ 0 & 1.5 \end{bmatrix} & \begin{bmatrix}  0 \\ 1 \end{bmatrix} \\  \begin{bmatrix} -.0005 & 0.001 \end{bmatrix} & [-0.127] \\
\end{array} \right ] 
\end{align}
The underlying system of (\ref{eq.gv}) is reconstructed with the random disturbances for inputs.  Note that $\mathbf{G}_{11}$ and $ \mathbf{G}_{12}$ are accurately reconstructed from the data.  The restricted Koopman operator also attempts to fit  $\mathbf{G}_{21}$ and $\mathbf{G}_{22}$ as a propagator on the random inputs, which will not be accurate by construction.  

If the controller has state-feedback, for example $u = -K x_2$ where $K = 1$.  The data in the last row becomes correlated with the second row.   In order to disambiguate the control from the $y_2$, a small disturbance is added to the input $u$ to only the snapshot matrix $\Omega$.  This provides the following approximate restricted Koopman operator:
\begin{align}
\mathbf{G} &=\left[ \begin{array}{cc} \mathbf{G}_{11}& \mathbf{G}_{12} \\  \mathbf{G}_{21}& \mathbf{G}_{22}  \end{array}\right] \approx
 \left [ \begin{array}{cc} \begin{bmatrix}0.1  & 0  \\ 0 & 1.5 \end{bmatrix} & \begin{bmatrix}  0 \\ 1 \end{bmatrix} \\  \begin{bmatrix} 0 & -1.5 \end{bmatrix} & [-1] \\
\end{array} \right ] 
\end{align}
where the dynamics on the controller now mimic the actual dynamics of $x_2$.  In this example, the restricted Koopman operator recovers the unstable underlying dynamics and discover that the inputs are being generated by a controller that is dependent on $x_2$.  

Consider the final input type:  the input has dynamics, but is not state dependent, for example $\dot{u} = -r u$ with $r = 0.01$ and $u(0) = 1$.  Similar to the other input types, we collect the data and find a restricted Koopman operator for the discrete system:
\begin{align}
\mathbf{G} &=\left[ \begin{array}{cc} \mathbf{G}_{11}& \mathbf{G}_{12} \\  \mathbf{G}_{21}& \mathbf{G}_{22}  \end{array}\right] \approx
 \left [ \begin{array}{cc} \begin{bmatrix}0.1  & 0  \\ 0 & 1.5 \end{bmatrix} & \begin{bmatrix}  0 \\ 1 \end{bmatrix} \\  \begin{bmatrix} 0 & 0 \end{bmatrix} & [0.99] \\
\end{array} \right ] .
\end{align}
The KIC architecture discovers the underlying dynamics of $\bx$ and the impact of $u$, but also finds the dynamics on $u$.  This perspective could be beneficial when considering multi-scale modeling where one scale is considered a forcing on another.  

The restricted KIC operator can be recovered from the data despite the unstable eigenvalue and various types of inputs.  Note that both the operator $\bA$ and $\bB$ are recovered from the underlying dynamical system (\ref{eq.gv}). The left Koopman modes, as in (\ref{eq.evProbt2}) of this operator are
\begin{align}
\mathbf{q} &=\left[ \begin{array}{cc} 1& 0 \\ 0 & 1 \end{array}\right] 
\end{align}
where these Koopman modes can be used to construct the eigenfunctions $\Psi_j$, described in \S\ref{ss:KICgen}.  A similar procedure can be utilized to find the right Koopman modes $\bv_j$ and eigenfunctions $\varphi_j$.  The right Koopman modes span both the states and the inputs.    
\subsection{Example 2 -- Nonlinear system with inputs}\label{sec:results:nonlinearsys}
We investigate how Koopman with control can be used to solve a nonlinear example with inputs.  In this example, we take the KIC form with the definition $* = \mathbf{0}$.  Consider the following nonlinear dynamical system from~\cite{Tu:2014a} and \cite{Brunton:2015Koopm}, but modified to include an input $u$
\begin{align}
\left [ \begin{array}{c} \dot{x}_1 \\ \dot{x}_2 \end{array} \right ] = 
 \left [ \begin{array}{c} \mu x_1 \\ \lambda (x_2 - x_1^2) + \delta u \end{array} \right ].
\end{align}
where $\lambda = 0.5$, $\mu = 2$, and $\gamma = 2$.  We use this example to investigate the effect of inputs or control on the nonlinear system.  The observable functions are carefully chosen, as in \cite{Brunton:2015Koopm}, to investigate this dynamical system given by the following:
\begin{align}
\left [ \begin{array}{c} y_1 \\ y_2 \\ y_3 \\ \tilde{u_1}\end{array} \right ] = 
\left [ \begin{array}{c} x_1 \\ x_2 \\ x_1^2 \\ u \end{array} \right ] ,~~~~~~
\left [ \begin{array}{c} \dot{y}_1 \\ \dot{y}_2 \\ \dot{y}_3 \end{array} \right ] =
 \left [ \begin{array}{ccc} \mu & 0 & 0   \\ 0 & \lambda & -\lambda  \\ 0 & 0 & 2 \mu \end{array} \right ] \left [ \begin{array}{c} y_1 \\ y_2 \\ y_3 \end{array} \right ]  + 
 \left [ \begin{array}{c} 0 \\ \delta \\ 0 \end{array} \right ] \tilde{u_1},
\end{align}
where the nonlinear function $y_3 = x_1^2$ has a convenient derivative which allows for closure of the dynamical system defined for the observables, see \cite{Brunton:2015Koopm} for more about closure of these dynamical systems.  We can transform the problem to include the inputs  
\begin{align}
\left [ \begin{array}{c} \dot{y}_1 \\ \dot{y}_2 \\ \dot{y}_3 \end{array} \right ] =
 \left [ \begin{array}{cccc} \mu & 0 & 0  & 0 \\ 0 & \lambda & -\lambda & \delta \\ 0 & 0 & 2 \mu & 0 \end{array} \right ] 
 \left [ \begin{array}{c} y_1 \\ y_2 \\ y_3 \\ \tilde{u_1} \end{array} \right ]  .
\end{align}
Now, we can collect measurement data in terms of the input and output variables $y_i$ and $u$.  In this case we used fifteen iterations with an initial conditions of $[5 2]^T$ . The restricted Koopman operator on these observables can be reconstructed:
\begin{align}
\mathbf{G} &=\left[ \begin{array}{cc} \mathbf{G}_{11}& \mathbf{G}_{12} \end{array}\right] \approx \left [ \begin{array}{ccccc}  2  & 0 & 0 & 0 \\ 0 & 0.5 & -0.5 & 2 \\ 0 & 0 & 4 & 0
\end{array} \right ] .
\end{align}
 The left Koopman modes can be constructed similar to \cite{Brunton:2015Koopm} and described in (\ref{eq.evProbt2}).  These Koopman modes $\bq_j$ can then be used to construct eigenfunctions $\Psi_j(\bx) = \langle \bx, \bq_j \rangle $.  These eigenfunctions span the Koopman operator for this nonlinear dynamical system.  The right Koopman modes and eigenfunctions can also be computed as described by (\ref{eq.evProbt2}).  Despite the nonlinear dynamical system, the KIC perspective constructs a linear dynamical system on the measurements  that can be used for prediction and control.
 
\subsection{Example 3 --  A biologically inspired nonlinear example
}\label{sec:results:SIR}
We investigate KIC on the classic Susceptible-Infected-Recovered (SIR) model.  This example contains a nonlinearity which is fundamentally different than for Example 2.  The nonlinearity does not have the same closure property.  Consider one version of the SIR models with inputs (represented by the vaccination of susceptible individuals):
\begin{align}
\left [ \begin{array}{c} \dot{S} \\ \dot{I} \\ \dot{R} \end{array} \right ] = 
 \left [ \begin{array}{c} -\beta S I +\nu (S + I + R) - \mu S - \text{Vacc}\\ 
 \beta S I - \gamma I - \mu I \\
 \gamma I - \mu R + \text{Vacc} \end{array}  \right ] 
\end{align}
\begin{figure}
\begin{center}
\begin{tabular}{cc}
\begin{overpic}[width=.425\textwidth]{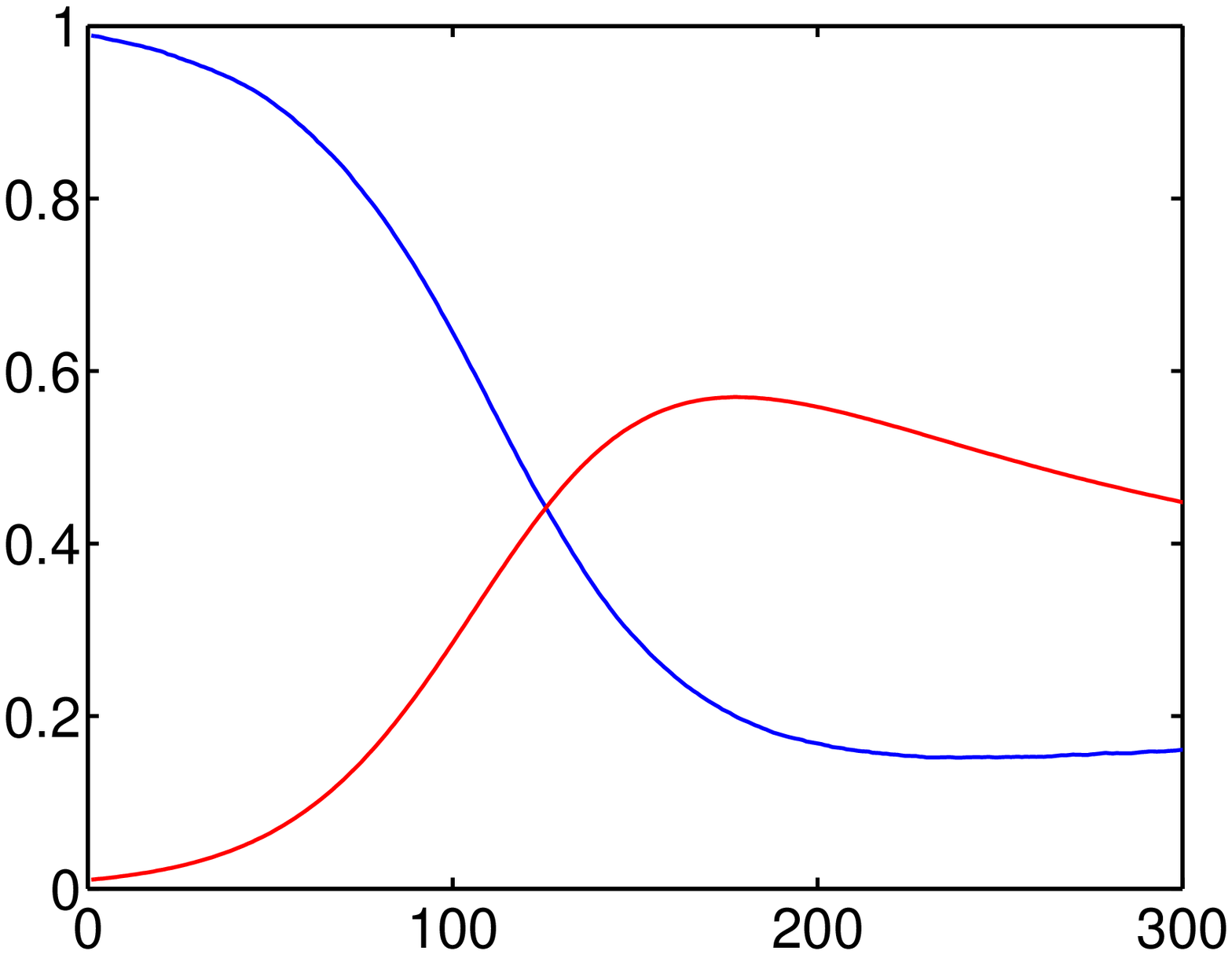}
\put(45,0){Time}
\put(45,60){SIR System}
\put(0,25){\rotatebox{90}{Population}}
\end{overpic}
&
\begin{overpic}[width=.425\textwidth]{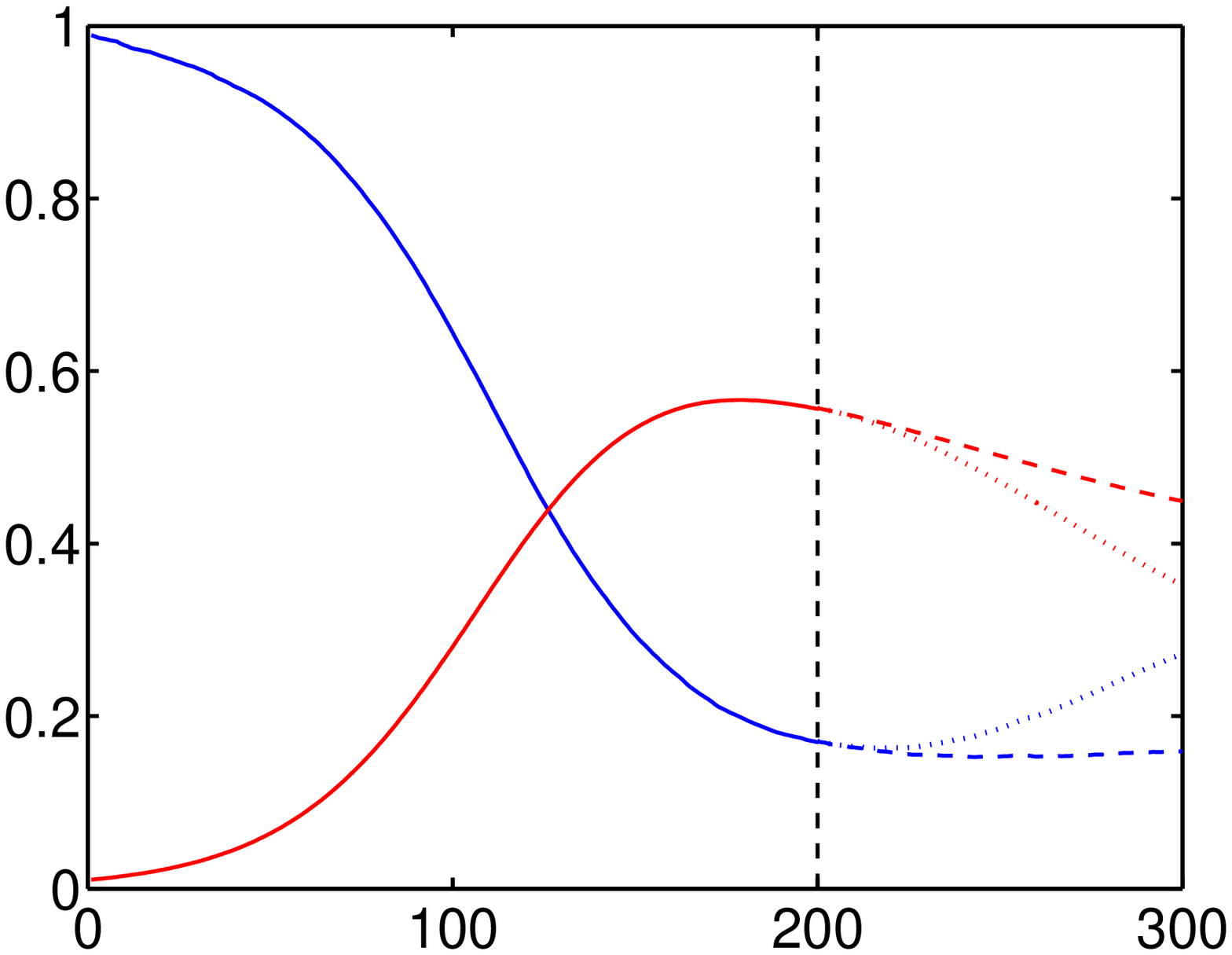}
\put(45,0){Time}
\put(45,60){Actual}
\put(66,60){Prediction}

\end{overpic}
\end{tabular}
\end{center}
\vspace{-.1in}
\caption{\small Left panel:  SIR dynamics with a 1\% seeding of infection at time zero.  Right panel:  The same SIR dynamics left of the dark dashed line.  To the right of the line,the dashed line indicates the KIC prediction with only the linear measurements $(S,I,R)$ as the output.  The dotted line indicates the prediction if the measurements $(S,I,R,SI)$ are used for the output.}\label{fig:SIR}
\end{figure}
where $\beta = 10$ is an infectious parameter, $\nu =1$ is a birthrate parameter depending on the total population of the community $S + I + R = 1$, $\mu = 1$ is the death rate, $\gamma =1$ is a recovery rate from infection, and $\text{Vacc}$ is a rate of vaccination.  The left panel shows the output of Fig.~\ref{fig:SIR} with seeding a 1\% infection at time zero and adding a small random amount of vaccination at each time step.  The nonlinearity in this example $SI$ is a mixed state quadratic nonlinearity.  We transform this continuous nonlinear dynamical system into a discrete linear dynamical system with a simple forward-euler scheme, augment the input space to include the nonlinearity $SI$ and inputs $\text{Vacc}$:
\begin{subequations}
\begin{align}
\text{Input:}~~~~~\left [ \begin{array}{ccccc} y_1  & y_2 & y_3 &y_4 & \gamma_1 \end{array} \right ]^T &= 
\left [ \begin{array}{ccccc} S & I & R & SI & \text{Vacc} \end{array} \right ]^T \label{eq.obi}
\end{align}
\begin{align}
\text{Output:}~~~~~~~~~~~~~~~\left [ \begin{array}{ccc} y_1  & y_2 & y_3  \end{array} \right ]^T &= 
\left [ \begin{array}{ccc} S & I & R  \end{array} \right ]^T 
\end{align}
\label{eq.obo}
\end{subequations}
giving the dynamical system:
\begin{subequations}
\begin{align}
\left [ \begin{array}{c} y_1 \\ y_2 \\ y_3 \end{array} \right ]_{k+1} &=
 \left [ \begin{array}{ccccc} 1/(\Delta t) - \mu +\nu & \nu & \nu &  -\beta & -1 \\ 
 0 & 1/(\Delta t)-\mu - \gamma & 0 & \beta & 0  \\ 
 0 & \gamma & 1/(\Delta t) - \mu & 0 & 1 
  \end{array} \right ] \left [ \begin{array}{c} y_1 \\ y_2 \\ y_3 \\ y_4 \\ \gamma_1 \end{array} \right ]_k   \label{eq:SIRd2}
  \end{align}
  \begin{align}
  \mathbf{Y}_o &= \mathbf{K} \mathbf{Y}_i,
  \label{eq:SIRd}
\end{align}
\end{subequations}
where $\mathbf{K}$ is the KIC operator from data, $\mathbf{Y}_o$ is the data in the output observables, and $\mathbf{Y}_i$ is the data in the input observables.  If the term $SI$ is included in the output observables, the derivative does not lend itself to a closed form.  The $\dot{SI} =\dot{S}I + S\dot{I} $ introduces the need for even more nonlinearities thus increasing the number of needed augmented observable functions, requiring $S^2I$, $SI^2$, and $I^2$.  We have not included a row in (\ref{eq:SIRd2}) for the time evolution of $y_4$ for this reason.  This introduces a practical difficulty in the implementation of this method on realistic complex systems.  For example, if we try to solve \eqref{eq:SIRd} with the augmented row as in Example 2, we do not find the correct operator.  Indeed, the last row has no semblance of the correct values.  The right panel of Fig.~\ref{fig:SIR} shows how this formulation (dotted line) provides an incorrect operator and thus is not able to accurately predict in the future after being trained on 200 time snapshots.  The dashed line shows the correct prediction after solving \eqref{eq:SIRd} with the input and output observables defined as in \eqref{eq.obi} and \eqref{eq.obo}.  

This particular equation can readily be solved with enough snapshot data.  Thus, choosing the correct observable functions becomes of paramount importance both on the input and output space.  A similar sentiment is expressed in \cite{Williams2015jnls}, but without considering separate input and output spaces.  Further, recent work has shown a statistical framework for determining which nonlinearities to include by sparsely choosing from a large library of possible dynamical terms \cite{Brunton:2015Sparse}.    

\subsection{Example 4 -- A nonlinear example with periodic solutions}
\label{sec:results:PO}

Here, we consider a nonlinear example that contains periodic solutions such that $\bx_k = \bx_{k+m} ~~\forall ~~k$.  The same example was considered in \cite{Rowley:2009} to illustrate that a common method, the discrete Fourier transform, for analyzing a periodic solutions can described in terms of the Koopman operator theory, more generally discussed \cite{Mezic:2005}. They illustrate that the Fourier expansion on the periodic orbit illustrates eigenfunctions of the Koopman operator.  

In this subsection, we consider a slightly different problem where the periodic orbit also has external inputs $\bu_k$, where each of the states $\bx_k$ and $\bu_k$ are contained in a set $S$.   Similar to \cite{Rowley:2009}, we can define a set of transformed states using the Fourier decomposition such that
\begin{align}
\bx_k = \sum_{j=0}^{m-1} \exp^{2\pi j k / m}\hat{\bx}_j,~~~k = 0,\dots,m-1.
\end{align}
Now, we consider expanding our input observables with a related transform called the Z-transform 
\begin{align}
\bu_k = \sum_{j=0}^{m-1} z^j \hat{\bu}_j,~~~k = 0,\dots,m-1.
\end{align}
where the Z-transform is a generalization of the Fourier decomposition such that in the infinite expansion $\sum_{k=-\infty}^\infty x_k z^{-k} = \sum_{k=-\infty}^\infty x_k \exp^{-jwn}$ for values of $|z| = 1$.  The Z-transform allows for a variety of different inputs and is significant when considering input-output transfer functions for linear systems.  Define a set of functions $\varphi_j(\bx_k,\bu_k) : S \rightarrow \mathbb{C}$ by 
\begin{align}
\varphi_j(\bx_k,\bu_k) = z^k, ~~~j,k = 0,\dots, m-1
\end{align}
then $\varphi_j$ are the right eigenfunctions of the Koopman operator $\bK$, with singular values $z$ on the unit circle and the left eigenfunctions $\psi_j$ will be defined by the Fourier decomposition 
\begin{align}
\bK \varphi_j(\bx_k,\bu_k) = \psi_j(f(\bx_k,\bu_k),\bu_{k+1}) = \psi_j(\bx_{k+1},\bu_{k+1}) = \exp^{2 \pi i j / m}\psi_j(\bx_k,\bu_k)
\end{align} 
Thus, the expansions of the input and output spaces can be written in terms of two different (but related) Hilbert spaces.  By restricting the phase space to just the periodic orbit $S$, the Koopman modes are given by the Fourier transform and the Z-transform.  Without inputs, the analysis reduces to the discrete Fourier transform.  With inputs, the choice of the z-transform allows for exogenous inputs to be considered in analyzing the periodic orbit.

\section{Discussion}\label{sec:discussion}

A wealth of modern applications are nonlinear and high-dimensional including distribution systems, internet traffic, and vaccinating of human populations in the developing world.  The need to develop quantitative and automatic methods to characterize and control these systems are of paramount importance to solving these large-scale problems.  In order to construct effective controllers, the complex system has to be well-understood.  In the case that we do not have well-established, physics-based governing equations, equation-free methods can help characterize these systems and offer insight into their control. 

Koopman operator theory and DMD offer a data-driven method to characterizing complex systems \cite{Mezic:2005,Rowley:2009}.   These methods are strongly grounded in the analysis of nonlinear systems and have been successfully applied in a number of fields such as fluid dynamics \cite{Mezic:2005,Schmid2008APS,Schmid:2009,Bagheri:2013}, epidemiology \cite{Proctor:2015EP}, video processing \cite{Grosek:2013}, and neuroscience \cite{Brunton2016}.  Further, this architecture has allowed for the incorporation of recent innovations from compressive sensing allowing insight into optimally measuring a system~\cite{Javanovic:2012,Tu:2014a,brunton:2014b}.  Generalizing Koopman for input-output systems allows for a broader set of  systems to be considered.  KIC is well-connected to DMDc, which is already having an impact analyzing input-output characteristics for systems with linear observables \cite{Proctor:2015DMDc,Dawson2015}.  The extension to Koopman operator theory allows for a larger set of observable functions to be included allowing for nonlinear system identification and the design of controllers.    

Theoretical innovations such as KIC will play an ever-increasing role in the characterization and control of complex systems.  We believe KIC and DMDc are well poised to be integrated in to a diverse set of engineering and science applications.  KIC is positioned to have a significant impact in the analysis and control of large-scale, complex, input-output systems.       

\section*{Acknowledgements} 

The authors would like to thank Bill and Melinda Gates for their active support of the Institute for Disease Modeling and their sponsorship through the Global Good Fund.  S. L. Brunton acknowledges support from the UW Mechanical Engineering department and from the eScience institute as a data science fellow.  J. N. Kutz acknowledges support from the U.S. Air Force Office of Scientific Research (FA9550-09-0174).  

\newpage
\bibliographystyle{plain}
\bibliography{DMDc}

\end{document}